\documentclass[12pt,a4paper]{article}
\usepackage{a4,amsmath,amssymb,amsfonts,amsthm}
\usepackage{amsrefs}
\usepackage[all]{xy}
 \usepackage[small,nohug,heads=vee]{diagrams}

\begin{document}
\bibliographystyle{plain}
 

\def\M#1{{\mathbb #1}}
\def\mR{\M{R}}           
\def\mZ{\M{Z}}           
\def\mN{\M{N}}           
\def\mQ{\M{Q}}       
\def\mC{\M{C}}  
\def\mG{\M{G}}
\def\mP{\M{P}}



\gdef\Spec{{\rm Spec}}
\def\rg{{\rm rg}}
\def\Hom{{\rm Hom}}
\def\Aut{{\rm Aut}}
 \def\Tr{{\rm Tr}}
 \def\Exp{{\rm Exp}}
 \def\Gal{{\rm Gal}}
 \def\End{{\rm End}}
 \def\det{{{\rm det}}}
 \def\Td{{\rm Td}}
 \def\ch{{\rm ch}}
 \def\che{{\rm ch}_{\rm eq}}
\def\Id{{\rm Id}}
\def\Zar{{\rm Zar}}
\def\Supp{{\rm Supp}}
\def\eq{{\rm eq}}
\def\Ann{{\rm Ann}}
\def\LT{{\rm LT}}
\def\Pic{{\rm Pic}}
\def\rg{{\rm rg}}
\def\et{{\rm et}}
\def\sep{{\rm sep}}
\def\ppcm{{\rm ppcm}}
\def\ord{{\rm ord}}
\def\Gr{{\rm Gr}}
\def\gr{{\rm gr}}
\def\ker{{\rm ker}}
\def\rk{{\rm rk}}
\def\Stab{{\rm Stab}}
\def\im{{\rm im}}
\def\Sm{{\rm Sm}}
\def\red{{\rm red}}
\def\Frob{{\rm Frob}}
\def\Ver{{\rm Ver}}


\def\beginProof{\par{\bf Proof.}}
 \def\endProof{${\qed}$\par\smallskip}
 \def\pr{^{\prime}}
 \def\prpr{^{\prime\prime}}
 \def\mtr#1{\overline{#1}}
 \def\ra{\rightarrow}
 \def\mfp{{\mathfrak p}}
 \def\mfm{{\mathfrak m}}
 
 \def\mQ{{\mathbb Q}}
 \def\mR{{\mathbb R}}
 \def\mZ{{\mathbb Z}}
 \def\mC{{\mathbb C}}
 \def\mN{{\mathbb N}}
 \def\mF{{\mathbb F}}
 \def\mA{{\mathbb A}}
  \def\mG{{\mathbb G}}
 \def\CI{{\cal I}}
 \def\CA{{\cal A}}
 \def\CE{{\cal E}}
 \def\CJ{{\cal J}}
 \def\CH{{\cal H}}
 \def\CO{{\cal O}}
 \def\CA{{\cal A}}
 \def\CB{{\cal B}}
 \def\CC{{\cal C}}
 \def\CC{{\cal C}}
 \def\CF{{\cal F}}
 \def\CL{{\cal L}}
 \def\CI{{\cal I}}
 \def\CM{{\cal M}}
\def\CP{{\cal P}}
 \def\CZ{{\cal Z}}
\def\CR{{\cal R}}
\def\CG{{\cal G}}
\def\CX{{\cal X}}
\def\CY{{\cal Y}}
\def\CV{{\cal V}}
\def\CW{{\cal W}}
 \def\wt#1{{\widetilde{#1}}}
 \def\mod{{\rm mod\ }}
 \def\refeq#1{(\ref{#1})}
 \def\blb{{\big(}}
 \def\brb{{\big)}}
\def\mc{{{\mathfrak c}}}
\def\mcpr{{{\mathfrak c}'}}
\def\mcprpr{{{\mathfrak c}''}}
\def\ss{{\rm ss}}
\def\parf{{\rm parf}}
\def\P1{{{\bf P}^1}}
\def\cod{{\rm cod}}
\def\pr{\prime}
\def\prpr{\prime\prime}
\def\ss{\scriptstyle}
\def\OX{{ {\cal O}_X}}
\def\mpartial{{\mtr{\partial}}}
\def\inv{{\rm inv}}
\def\indlim{\underrightarrow{\lim}}
\def\prolim{\underleftarrow{\lim}}
\def\pprolim{'\prolim'}
\def\Pro{{\rm Pro}}
\def\Ind{{\rm Ind}}
\def\Ens{{\rm Ens}}
\def\without{\backslash}
\def\pbdb{{\Pro_b\ D^-_c}}
\def\qc{{\rm qc}}
\def\Com{{\rm Com}}
\def\an{{\rm an}}
\def\gfield{{\rm\bf k}}
\def\s{{\rm s}}
\def\dR{{\rm dR}}
\def\ari#1{\widehat{#1}}
\def\ul#1{\underline{#1}}
\def\sul#1{\underline{\scriptsize #1}}
\def\mou{{\mathfrak u}}
\def\ich{\mathfrak{ch}}
\def\cl{{\rm cl}}
\def\K{{\rm K}}
\def\R{{\rm R}}
\def\F{{\rm F}}
\def\L{{\rm L}}
\def\pgcd{{\rm pgcd}}
\def\rc{{\rm c}}
\def\N{{\rm N}}
\def\E{{\rm E}}
\def\H{{\rm H}}
\def\CHOW{{\rm CH}}
\def\A{{\rm A}}
\def\d{{\rm d}}
\def\Res{{\rm  Res}}
\def\GL{{\rm GL}}
\def\Alb{{\rm Alb}}
\def\alb{{\rm alb}}
\def\Hdg{{\rm Hdg}}
\def\Num{{\rm Num}}
\def\Irr{{\rm Irr}}
\def\Frac{{\rm Frac}}
\def\Sym{{\rm Sym}}
\def\TV{\rm TV}
\def\indlim{\underrightarrow{\lim}}
\def\prolim{\underleftarrow{\lim}}
\def\Ver{{\rm Ver}}
\def\hn{{\rm hn}}
\def\min{{\rm min}}
\def\max{{\rm max}}
\def\Div{{\rm Div}}
\def\sm{{\rm sm}}
\def\Mor{{\rm Mor}}
\def\wwt#1{\wt{\wt{#1}}}
\def\h{{\rm h}}
\def\charac{{\rm char}}


\def\RHom{{\rm RHom}}
\def\rRHom{{\mathcal RHom}}
\def\rHom{{\mathcal Hom}}
\def\dotimes{{\overline{\otimes}}} 
\def\Ext{{\rm Ext}}
\def\rExt{{\mathcal Ext}}
\def\Tor{{\rm Tor}}
\def\rTor{{\mathcal Tor}}
\def\SP{{\mathfrak S}}
\def\perf{{\rm perf}}
\def\T{{\rm T}}
\def\H{{\rm H}}
\def\D{{\rm D}}
\def\Del{{\mathfrak D}}
\def\II{{\rm II}}
\def\sh{{\rm sh}}
 \def\length{{\rm length}}
 \def\min{{\rm min}}
 \def\max{{\rm max}}

 \newtheorem{theor}{Theorem}[section]
 \newtheorem{prop}[theor]{Proposition}
 \newtheorem{propdef}[theor]{Proposition-Definition}
 \newtheorem{sublemma}[theor]{sublemma}
 \newtheorem{cor}[theor]{Corollary}
 \newtheorem{lemma}[theor]{Lemma}
 \newtheorem{sublem}[theor]{sub-lemma}
 \newtheorem{defin}[theor]{Definition}
 \newtheorem{conj}[theor]{Conjecture}
 \newtheorem{quest}[theor]{Question}

 \parindent=0pt
 \parskip=5pt

 \author{Henri GILLET\footnote{Department of Mathematics, Statistics, and Computer Science (m/c 249), 
University of Illinois at Chicago, 
851 South Morgan Street, Chicago, IL 60607-7045, USA}
 \, \, \&\ \ Damian R\"OSSLER\footnote{Mathematical Institute
University of Oxford, 
Andrew Wiles Building, 
Radcliffe Observatory Quarter, 
Woodstock Road, 
Oxford, 
OX2 6GG, 
United Kingdom}}
 \title{Rational points of varieties with ample cotangent bundle over function fields}
\maketitle

\begin{abstract}
Let $K$ be the function field of a smooth curve over an algebraically closed field $k$. Let $X$ be a scheme, which is smooth and projective over $K$. Suppose that the cotangent bundle $\Omega_{X/K}$ is ample. Let \mbox{$R:=\Zar(X(K)\cap X)$} be the Zariski closure of the set of all $K$-rational points 
of $X$, endowed with its reduced induced structure. We prove that for each irreducible component ${\mathfrak R}$ of $R$, 
 there is a projective variety ${\mathfrak R}'_0$ over $k$ and a finite  and surjective 
$K^\sep$-morphism ${\mathfrak R}'_{0,K^\sep}\to {\mathfrak R}_{K^\sep}$, which is birational when 
$\charac(K)=0$. 

\noindent This improves on results of Noguchi and Martin-Deschamps in characteristic $0$. 
In positive characteristic, our result can be used to give the first examples of varieties, which are not embeddable in abelian varieties and satisfy an analog of the Bombieri-Lang conjecture. 
\end{abstract}

\section{Introduction}

Recall that the Bombieri-Lang conjecture (see \cite[middle of p. 108]{Noguchi-A-higher})  asserts that the set of rational points of a 
variety of general type over a number field is not dense. This 
conjecture can be proven in the situation where the variety is embeddable 
in an abelian variety (this is a deep result of Faltings, see \cite{Faltings-Diophantine}) but it is not known to be true in any  
other situation, as far as the authors know. 

Over function fields, it seems reasonable to make the following conjecture, which 
must have been part of the folklore for some time.

\begin{conj}['Bombieri-Lang' conjecture over function fields]
Let $K_0$ be the function field of a variety over an algebraically closed field 
$k_0$. Let $Z$ be a variety of general type over $K_0$. Suppose that $\Zar(Z(K_0))=Z$. Then there exists a variety $Z_0$ over $k_0$ 
and a rational, dominant, generically finite $K_0^\sep$-map 
$g:Z_{0,K^\sep_0}\to Z_{K^\sep_0}.$
\label{BLpc}
\end{conj}
One might speculate that the rational map $g$ appearing in Conjecture \ref{BLpc} can also be taken to be generically purely inseparable. 

In his article \cite[p. 781]{Lang-Higher}, Lang gave a loose formulation of 
Conjecture \ref{BLpc} for 
$\charac(k_0)=0$. 
If $Z$ is embeddable in an abelian variety and $\charac(k_0)=0$, Conjecture \ref{BLpc}  can be proven (see \cite{Buium-Intersections} and 
\cite{FW-Rational}). It can also be proven in the situation where $\charac(k_0)=0$, the variety is smooth and its cotangent bundle is ample. This 
is a result of Noguchi, which was also proved independently by Martin-Deschamps (see \cite{Noguchi-A-higher} and \cite{Martin-Deschamps-Proprietes}). 
When $Z$ is embeddable in an abelian variety and $\charac(k_0)>0$, 
Conjecture \ref{BLpc} is a theorem of Hrushovski (see \cite{Hrushovski-Mordell-Lang}). 
See also \cite{Rossler-MMML}, \cite{Ziegler-Mordell-Lang} and \cite{Benoist-Bouscaren-Pillay-MLMM} for 
different proofs of Hrushovski's theorem. When $Z$ is of dimension $1$ and $\charac(k_0)=0$, Conjecture 
\ref{BLpc} was first proven by Manin and Grauert and several 
other proofs were given in the course of the 1970s (eg by Parshin and Arakelov). See the articles 
\cite{Grauert-Mordell}, \cite{Manin-Rational}, \cite{Manin-Letter} and \cite[chap. I]{Cornell-Silverman-Arithmetic}. When $Z$ is of dimension $1$ and $\charac(k_0)>0$, 
Conjecture \ref{BLpc} was first proven by Samuel and other proofs were given 
later by Szpiro and Voloch. See the articles \cite{Samuel-Complements}, \cite{Szpiro-Seminaire-Pinceaux} and 
\cite{Voloch-Towards}. 
 
In the following paper, we shall prove Conjecture \ref{BLpc} in the situation where $K_0$ has transcendence degree $1$ over its prime field and $Z$ is a subvariety of a larger variety $Z'$, where $Z'$ is smooth and has ample cotangent bundle over over $K_0$. See Theorem \ref{Mtheor} below. Theorem \ref{Mtheor} can be used to give the first examples (to the authors knowledge)  of 
varieties of general type in positive characteristic that satisfy Conjecture \ref{BLpc} and are not embeddable in abelian varieties (see \cite{Dupuy-Examples} 
about this). Theorem \ref{Mtheor} also provides a strengthening of Noguchi's result in the situation where $K_0$ has transcendence degree $1$ over $\mQ$ 
(because unlike Noguchi we do not assume that $Z=Z'$).

Here is a precise formulation of our result.

Let $K$ be the function field of a smooth curve $U$ over an algebraically closed field 
$k$. Let $X$ be a scheme, which is smooth and projective over $K$. 

We prove: 

\begin{theor} Suppose that the cotangent bundle $\Omega_X:=\Omega_{X/K}$ is ample.  
Let \mbox{$R:=\Zar(X(K)\cap X)$} be the Zariski closure of the set of all $K$-rational points 
of $X$, endowed with its induced reduced structure. For each irreducible component ${\mathfrak R}$ of $R$, 
there is a projective variety ${\mathfrak R}'_0$ over $k$ and a finite and surjective 
$K^\sep$-morphism $h:{\mathfrak R}'_{0,K^\sep}\to {\mathfrak R}_{K^\sep}$. 

If $\charac(k)=0$ then there exists a morphism 
$h$ as above, which is birational. 

If $R={\mathfrak R}=X$ and 
$$
\charac(k)>\dim(X)^2\int_X{\rm c}_1(\Omega_X)^{\dim(X)}
$$
then there exists a morphism $h$ as above, which is birational (and thus an isomorphism).

\label{Mtheor}
\end{theor}
Here ${\rm c}_1(\Omega_X)$ denotes the first Chern class of the vector bundle $\Omega_X$ in 
the Chow intersection ring ${\rm CH}^\bullet(X)$ of $X$ and 
$$\int_X:{\rm CH}^{\dim(X)}(X)\to{\rm CH}^0(\Spec\, K)\simeq\mZ$$ refers to 
the push-forward morphism. 

Note that the number $\dim(X)^2\int_X{\rm c}_1(\Omega_X)^{\dim(X)}$ is always positive, 
since $\Omega_X$ is ample.

{\bf Remark.} In Theorem \ref{Mtheor}, we assume that 
$K$ is the function field of a curve (see also the discussion above). This restriction, which can probably be removed at the price 
of added technicality, 
comes from the fact that we need to consider a smooth compactification of the curve $U$ in the proof and 
also from the fact that our proof depends on the existence of N\'eron desingularisations.  It seems difficult to 
reduce the case of a high dimensional base to the case of a base of dimension one by a slicing argument, 
because it is not clear (to us) how the assumption that the rational points are Zariski dense in a certain closed 
set behaves when one passes to a fibre in a family.

We now describe the strategy of the proof of Theorem \ref{Mtheor}, which can be viewed 
as a refinement of the method of Grauert (see \cite[chap. VI]{Lang-Survey} for a nice overview of this method), in which all the higher jet schemes 
are brought into the picture (unlike Grauert, who considers only the first jet scheme).

Here is how Grauert and his followers proceed in characteristic $0$, 
in the situation where the rational points are dense in the whole variety.  One first shows that 
the rational points of $X$  lift to rational points of 
the first jet scheme. Next, one shows that the ``non-constant'' rational points concentrate on a closed subscheme $\Sigma$ (say), which is finite and generically inseparable over $X$. To establish this last fact, one needs to consider a non-singular 
compactification of $X/K$ over a compactification of $U$ and use the height machine over function fields.
The proof can now be completed quickly because 
generically inseparable morphisms are birational in characteristic $0$ and 
thus the projection of $\Sigma$ onto $X$ is an isomorphism. 
In view of the definition of the first jet scheme, this means that the Kodaira-Spencer 
class of $X$ vanishes. Using the exponential map, one can conclude 
from this  that $X$ descends to $k$ (up to a separable 
extension of $K$).

If one tries to carry through the above proof in positive characteristic, the first problem that one faces 
 is that one cannot easily construct a non-singular compactification of $X$, unless one 
assumes the existence of resolutions of singularities in positive characteristic. 

Next, even if one supposes that this first problem can be solved, one is faced with the basic problem that 
the projection $\Sigma\to X$ might not be an isomorphism. 

Finally, even if the projection $\Sigma\to X$ can be shown to be an isomorphism, ie even if the Kodaira-Spencer class 
of $X$ vanishes, one cannot conclude that $X$ descends to $k$. For example, 
any smooth proper curve over $K$, which descends to $K^p$, has a vanishing Kodaira-Spencer 
class. 

Here is how we deal with these issues. For the first problem, we replace the non-singular compactification 
by a N\'eron desingularisation, which is not compact, but suffices for our purposes. The second and third issues are dealt with 
simultaneously. We show that after a finite purely inseparable base-change the entire tower 
of jet schemes becomes trivial, in the following sense: the base-change of the first jet scheme 
has a section, the base-change of the second jet scheme has a section over the image of the first 
section, the base-change of the third jet scheme has a section over the image of the second section and 
so on. We show that the trivialisation of the tower of jet schemes is a consequence of a generalisation of a cohomological result of Szpiro and Lewin-M\'en\'egaux (see before Proposition \ref{propn} below). Our proof of this generalisation is not based on the same 
principle as the result of Szpiro and Lewin-M\'en\'egaux. Specializing all this to a closed point $u_0$ of $U$, we obtain a morphism of formal schemes between a constant formal scheme and the completion at $u_0$ of (a suitable model of) $X$. Applying Grothendieck's formal GAGA theorem and using 
the fact that the completion of $U$ at $u_0$ is an excellent discrete valuation ring, we can construct  
the required morphism ${\mathfrak R}'_{0,K^\sep}\to {\mathfrak R}_{K^\sep}.$ 
The condition given in the theorem for this morphism to be birational in positive characteristic 
is a consequence of a result of Langer (see below for detailed references) and of the fact that the minimal slope 
of $\Omega_X$ must be positive, since $\Omega_X$ is ample. 

The method that we just outlined also allows us to treat the situation where the rational 
points are not Zariski dense (ie when $R\not =X$). This was (apparently) not accessible before even in characteristic $0$. 
Note that when $R$ is smooth then $R$ will also have an ample cotangent bundle but we do not need 
to make this assumption here (the fact that the ambient variety $X$ is smooth with ample cotangent bundle suffices).

Here is the structure of the text.

In subsection \ref{ssec1}, we recall various facts about the geometry of torsors under vector bundles (in particular, ample 
vector bundles). In subsection \ref{ssec2}, we prove an injectivity criterion for purely inseparable pull-back maps between 
first cohomology groups of vector bundles (Corollary \ref{corinj}) and we prove a basic vanishing result (Proposition \ref{propn}) for the group of 
global sections of a coherent sheaf, which is twisted by a sufficiently high power of Frobenius pull-backs of an ample bundle. 
In section \ref{sec3}, we prove Theorem \ref{Mtheor}. As explained above, our proof does not use the exponential map 
but uses formal schemes directly and thus differs in nature from the proofs of Noguchi and Martin-Deschamps (see op. cit.) 
even when $\charac(k)=0$. 

We recommend to the reader to first read the proof with the supplementary assumption that $U$ is proper over $k$ and that $X$ extends to a smooth and projective scheme over $U$. Many technicalities of the proof disappear when that (unrealistic\dots) supplementary assumption is made. A more technical outline of the proof is given at the beginning of section \ref{sec3}.  

Finally, we would like to extend our heartfelt thanks to the anonymous referee for his thorough checking of the text. Without him, this article would be a lot less clearly written.

{\bf Notations.} A scheme of positive characteristic is a scheme $S$, such that 
for all points $s\in S$, the local ring $\CO_s$ is a ring of positive characteristic. A commutative ring 
$R$ is of positive characteristic if there exists a prime number $p_0$ such that $p_0\cdot 1_R=0$, where 
$1_R$ is the unit element of $R$. If $S$ is a scheme of positive characteristic, we write $F_S$ for 
the absolute Frobenius endomorphism of $S.$ The acronym wrog stands 
for ``without restriction of generality''. If $Y$ is an integral scheme, we write 
$\kappa(Y)$ for the function field of $Y$. ``Almost all'' means ``for all but a finite number''.

We are deeply grateful to the anonymous referee for his thorough checking of the text. Without him, this article would be a lot less clearly written.

\section{Preliminaries}

\label{sec2}

\subsection{The geometry of the compactifications of torsors under vector bundles}

\label{ssec1}
In this subsection, we recall various results proven in \cite[par. 1]{Martin-Deschamps-Proprietes}. 

Let $S$ be a scheme, which is  of finite type over a field $k_0$. 

If $V$ a locally free sheaf over $S$, we shall write $\mP(V)$ for the $S$-scheme representing 
the functor on $S$-schemes
$$
T\mapsto \{\textrm{\rm iso. classes of surjective morphisms of $\CO_T$-modules $V_T\to Q$,}
$$
$$
\textrm{where $Q$ is locally free of rank $1$}\}.
$$
By construction, $\mP(V)$ comes with a universal line bundle $\CO_P(1)$. 
Let now
\begin{equation}
\CE:0\to \CO_S\to E\to F\to 0
\label{ES}
\end{equation}
be an exact sequence of locally free sheaves over $S$. Consider the $S$-group scheme $\underbar F:=\Spec(\Sym(F))$ representing the group 
functor on $S$-schemes sending $T$ to $F^\vee_T(T)$. 
Let $R_\CE$ be the functor from $S$-schemes to sets given by 
\begin{equation}
T\mapsto\{\textrm{morphisms of $\CO_T$-modules $E_T\mapsto \CO_T$ splitting $\CE_T$}\}.
\label{UF}
\end{equation}
There is an obvious (group functor-)action of $\underbar F$ on $R_\CE$. 

\begin{itemize}
\item[(i)] (see \cite[Prop. 1 and proof]{Martin-Deschamps-Proprietes}) The natural morphism $\mP(F)\to \mP(E)$ is a closed immersion and there is an isomorphism of line bundles 
$\CO(\mP(F))\simeq\CO_P(1)$.
\item[(ii)]  (see \cite[Prop. 1]{Martin-Deschamps-Proprietes}) The complement $\mP(E)\backslash \mP(F)$ represents the functor 
$R_\CE$. The isomorphism of functors on $S$-schemes $R_\CE\to \mP(E)\backslash \mP(F)$ can be 
described as follows. There is a  natural transformation of functors $R_\CE\to \mP(E)$ sending a 
morphism of $\CO_T$-modules $E_T\mapsto \CO_T$ splitting $\CE_T$ to the same morphism $E_T\to\CO_T$, viewed 
as a morphism from $E_T$ onto a locally free sheaf of rank $1$ (the latter being the trivial sheaf). This gives a morphism of schemes  $R_\CE\to \mP(E)$, which is an open immersion onto $\mP(E)\backslash \mP(F)$. 
\end{itemize}
Thus
\begin{itemize}
\item[(iii)] the scheme $R_\CE$ with its $\underbar F$-action is an $S$-torsor under 
$\underbar F$ (see also \cite[top of p. 42]{Martin-Deschamps-Proprietes}). 
\end{itemize}
Further, by (i): 
\begin{itemize}
\item[(iv)] if $E$ is ample and $S$ is proper over $k_0$ then the scheme $\mP(E)\backslash \mP(F)$ is affine
\end{itemize}
(point (iv) will actually not be used in the text). 

Let us now suppose until the end of this section that 
$F$ is ample and that $S$ is proper over $k_0$.

\begin{itemize}
\item[(v)] (see \cite[Prop. 2]{Martin-Deschamps-Proprietes}) if $Z\hookrightarrow R_\CE$ is a subscheme, which is closed in $\mP(E)$, then the induced map $Z\to S$ is finite 
and has only a finite number of fibres that contain more than one point.
\item[(vi)] (see \cite[Prop. 3]{Martin-Deschamps-Proprietes}) for all sufficiently large $n\in\mN$, the line bundle $\CO_P(n+1)$ is generated by its global sections 
and in this case the induced $k_0$-morphism $$\phi_{n}:\mP(E)\to\mP(\Gamma(\CO_P(n+1)))\simeq\mP_{k_0}^n$$ is generically finite;
\item[(vii)] (see \cite[Prop. 3]{Martin-Deschamps-Proprietes}) if the line bundle $\CO_P(n)$ is generated by its global sections, the positive-dimensional fibres of the morphism $\phi_n$ are disjoint from $\mP(F)$ (where $\phi_n$ is as in (vi)).
\end{itemize}
From the fact that fibre dimension is upper semi-continuous (see \cite[IV, 13.1.5]{EGA}) and (vii) we deduce that 
\begin{itemize}
\item[(viii)] the union $I(\phi_n)$ of the positive dimensional fibres of $\phi_n$ is closed in $\mP(E)$ and 
is contained in $R_\CE$. 
\end{itemize}
We endow $I(\phi_n)$ with its reduced-induced structure. From (v) we deduce that 
\begin{itemize}
\item[(ix)] the morphism $\lambda_n:I(\phi_n)\to S$ is finite and has only a finite number of fibres that contain more than one point; in particular, if $T\subseteq S$ is irreducible and closed, then 
$I(\phi_n)\cap\lambda_n^{-1}(T)$ has at most one irreducible component  
dominating $T$, and if it exists, this component is generically radicial over $T$.
\end{itemize}

The following diagram summarises the relationships between and the various properties of the various morphisms 
introduced in this subsection.
\vskip0.5cm
\centerline{
\xymatrix{
\mP(F)\ar@{^{(}->}[r]\ar[dddr] & \mP(E)\ar[r]^{\phi_n \textrm{\,\,\,\,\,\,\,}}_{\textrm{gen. finite\,\,\,\,\,}} & \mP(\Gamma(\CO_P(n+1)))\ar@/^4pc/[ddddl]\\
& \,\,\,\,\,\,\,\,\,\,\,\,\,\,\,\,\,\,\,\,\,\,\,\,\,\,\,\,\,\, \mP(E)\without\mP(F)\ar[u]^{\textrm{open imm.}}\textrm{\tiny\,\, torsor under $\underbar F$} &\\
& I(\phi_n)\ar@{_{(}->}[u]\ar[d]^{\textrm{finite, gen. rad.}} &\\
& S\ar[d] & \\
& \Spec(k_0)&
}}
\vskip0.5cm

We shall also record the following geometric consequence of (v):

\begin{itemize}
\item[(x)] If $Z\hookrightarrow R_\CE$ and $Z'\hookrightarrow R_\CE$ are integral subschemes, such that the induced morphisms $Z\to S$ and $Z'\to S$ are finite and surjective then $Z=Z'$.
\end{itemize}
In particular, if $S$ is integral, then any two sections of $R_\CE$ must coincide. 

Finally, we shall also need the
\begin{itemize}
\item[(xi)] Every torsor under $\underbar F$ is isomorphic to a torsor $R_\CE$ for some 
exact sequence $\CE$ as in \refeq{ES}. The class in $H^1(S,F^\vee)\simeq \Ext^1(\CO_S,F^\vee)$ corresponding to $R_\CE$ is 
the image of $1\in H^0(S,\CO_S)$ in $H^1(S,F^\vee)$ under the connecting map in the 
long exact sequence
$$
0\to H^0(S,F^\vee)\to H^0(S,E^\vee)\to H^0(S,\CO_S)\to H^1(S,F^\vee)\to\dots
$$
associated with the dual exact sequence $\CE^\vee$. 
\end{itemize}

We leave the details of the proof of (xi) (which is but unravelling the definitions) to the reader. 
The assumptions that $S$ is proper over $k_0$ and that $F$ is ample are not necessary in (xi).

\subsection{Torsors under vector bundles and purely inseparable base-change}

\label{ssec2}

If $W$ is a quasi-coherent $\CO_Y$-module on a integral scheme $Y$, we shall write
$$
\Gamma(Y,W)_g:=\{e\in W_{\kappa(Y)}\,|\, 
\exists\, \sigma\in\Gamma(Y,W):\sigma_{\kappa(Y)}=e\}.
$$

\begin{lemma} 
Let $Y$ be a normal, noetherian and integral scheme. Let $W$ be a vector bundle over $Y$. 
Let $T\to Y$ be a torsor under $W$ and let $Z\hookrightarrow T$ be a closed immersion, 
where $Z$ is an integral scheme. 
Suppose that the induced morphism $f:Z\to Y$ is quasi-finite, radicial and dominant.  Suppose that $\Gamma(Z,\Omega_f\otimes f^*W)_g=0$. Then $f$ is an open immersion.
\label{allimp}
\end{lemma}
\beginProof \,
Let $\pi:T\times_Y T\to Y$. 
We consider the scheme $T\times_Y(T\times_Y T)$. Via the projection 
on the second factor $T\times_Y T$, this scheme is naturally a torsor under the vector bundle $\pi^*W$ . This torsor has two sections: 

- the section $\sigma_1$ defined by the formula $t_1\times t_2\mapsto t_1\times (t_1\times t_2)$;

- the section $\sigma_2$ defined by the formula $t_1\times t_2\mapsto t_2\times (t_1\times t_2)$.

Since $T\times_Y(T\times_Y T)$ is a torsor under $\pi^*W$, there is a section 
$s\in \Gamma(T\times_Y T,\pi^*W)$ such that $\sigma_1+s=\sigma_2$ and by construction 
$s(t_1\times t_2)=0$ iff $t_1=t_2$. In other words, $s$ vanishes precisely on  
the diagonal of $T\times_Y T$. 

Now let $\Delta_Z:Z\hookrightarrow Z\times_Y Z$ be the diagonal immersion. 
Let $\Delta_Z^{(1)}:Z^{(1)}\hookrightarrow Z\times_Y Z$ be the first infinitesimal 
neighborhood of $\Delta_Z$. By the definition of the differentials, we have an exact sequence of 
sheaves
$$
0\to\Omega_f\to \CO_{Z^{(1)}}\to \CO_Z\to 0
$$
which gives rise to an exact sequence
$$
0\to\Gamma(Z,\pi^*W|_Z\otimes \Omega_f)\to\Gamma(Z^{(1)},\pi^*W|_{Z^{(1)}})\to 
\Gamma(Z,\pi^*W|_Z)
$$
Now by construction, the image of the section $s\in \Gamma(T\times_Y T,\pi^*W)$  
in  $\Gamma(Z,\pi^*W|_Z)$ vanishes. Hence the image of the section $s\in \Gamma(T\times_Y T,\pi^*W)$  
in  $\Gamma(Z^{(1)},\pi^*W|_{Z^{(1)}})$ is the image of a section $s_0\in \Gamma(Z,\pi^*W|_Z\otimes \Omega_f)$. 
By assumption, we have $s_{0,\kappa(Z)}=0$ and thus by the construction of $s$ the immersion 
$Z\hookrightarrow Z^{(1)}$ is generically an isomorphism. Hence $\Omega_{f,\kappa(Z)}=0$ and thus 
the extension of function fields $\kappa(Z)|\kappa(Y)$ must be of degree 
one because it is a purely inseparable extension. Thus $f$ is birational. 
Zariski's main theorem now implies that $f$ is an open immersion.
\endProof

\begin{cor} Let $Y$ be a normal, noetherian and integral scheme of positive characteristic. Suppose 
that $F_Y$ is a finite morphism. Let $W$ be a vector bundle over $Y$. 
If $\Gamma(Y,\Omega_{F_Y}\otimes F^*_Y(W) )_g=0$ then the natural map of abelian groups 
$H^1(Y,W)\to H^1(Y,F_Y^*W)$ is injective.
\label{corinj}
\end{cor}
Recall that $F_Y:Y\to Y$ is the absolute Frobenius morphism (see the conventions on notations 
at the end of the introduction).
\beginProof\, 
Consider an element in the kernel of the map 
$H^1(Y,W)\to H^1(Y,F_Y^*W)$. Let $T\to Y$ be a torsor under $W$ corresponding 
to this element. By assumption, we have a diagram 

\centerline{
\xymatrix{
F_Y^*T\ar[r]^{J}\ar@<-.5ex>[d] & T\ar[d]\\
Y\ar[r]^{F_Y}\ar@<-.5ex>[u]_{\sigma} & Y\\
}
}
where the square is cartesian and $\sigma$ is a $Y$-morphism.
Now define $Z$ as the scheme-theoretic image $J_*(\sigma(Y))$ by $J$. 
By construction, we have a factorisation $Y\stackrel{\phi}{\to} Z\stackrel{f}{\to} Y$, where the arrow $f$ is the natural projection and $f\circ\phi=F_Y$. Now using \cite[Th. 26.5, p. 202]{Matsumura-Commutative}, we see that we have $$p^{\rk_{\kappa(Z)}(\Omega_{f,\kappa(Z)})}=[\kappa(Z):\kappa(Y)],$$
$$
p^{\rk_{\kappa(Y)}(\Omega_{F_Y,\kappa(Y)})}=[\kappa(Y)^{p^{-1}}:\kappa(Y)]
$$
and
$$
p^{\rk_{\kappa(Y)}(\Omega_{\phi,\kappa(Y)})}=[\kappa(Y):\kappa(Z)].
$$
Thus in the standard exact sequence
$$
\phi^*\Omega_{f,\kappa(Y)}\to \Omega_{F_Y,\kappa(Y)}\to\Omega_{\phi,\kappa(Y)}\to 0.
$$
the leftmost arrow is injective. 
Now from this and from the fact that \mbox{$\Gamma(Y,\Omega_{F_Y}\otimes F^*_Y(W) )_g=0,$} we  conclude that 
$\Gamma(Z,\Omega_{f}\otimes f^*W)_g=0$. Thus, by Lemma \ref{allimp}, $f$ is an isomorphism and thus the torsor 
$T\to Y$ is trivial.\endProof

Note that if $Y$ is projective over an algebraically closed field, then a weaker form of 
Corollary \ref{corinj} (which is not sufficient for our purposes) is contained in \cite[exp. 2, Prop. 1]{Szpiro-Seminaire-Pinceaux}. The proof 
given there depends on the existence of the Cartier isomorphism. 

Let now $S$ be an integral scheme, which is projective and geometrically reduced over a field $k_0$ of characteristic $p$. 
As before, let $F_S:S\to S$ be the absolute Frobenius morphism on $S$. 

\begin{prop}
Suppose that $\dim(S)>0$. Let $V$ be an ample vector bundle of rank $r$ over $S$. Let $W$ be a 
coherent sheaf over 
$S$.  Then we have 
\mbox{$H^0(S,F^{n,*}_S(V^\vee)\otimes W)_g=0$} for all sufficiently large $n\geqslant 0$.
\label{propn}
\end{prop}
\beginProof\, We may assume wrog that $k_0$ is algebraically closed. 

The proof is by induction on the dimension $d\geqslant 1$ of $S$.

Suppose that $d>1$. Consider a pencil of hypersurfaces in $S$ and let 
$b:\wt{S}\to S$ the total space of the pencil. By construction, we are given a morphism $m:\wt{S}\to\mP^1_{k_0}$. 
Let $\eta\in \mP^1_{k_0}$ be the generic point. Let $n_0$ be sufficiently large, so that 
$$H^0(\wt{S}_\eta,F^{n,*}_{\wt{S}_\eta}((b^*V^\vee)_\eta)\otimes (b^*W)_\eta)_g=0$$ for all $n\geqslant n_0$. This is possible by the induction hypothesis and 
because $(b^*V)_\eta$ is ample. The fact that $(b^*V)_\eta$ is ample is a consequence of the fact that the restriction of $V$ to any closed fibre of $m$ is ample and of the fact that 
ampleness on the fibre of $m$ is a constructible property (see \cite[IV, §9, Cor. 9.6.4]{EGA}). Now if we had $H^0(S,F^{n,*}_S(V^\vee)\otimes W)_g\not=0$ for 
some $n\geqslant n_0$ then 
the pull-back $b^*(F^{n,*}_S(V^\vee)\otimes W)$ would have a section that does not vanish 
at the generic point of $\wt{S}_\eta$, which is a contradiction.

Thus we are reduced to prove the statement for $d=1$. We may replace wrog 
$S$ by its normalisation. Since $S$ is now a non-singular curve, we know that $V$ is cohomologically $p$-ample (see \cite[Rem. 6), p. 91]{Migliorini-Some}). 
Also, $W$ is now the direct sum of a torsion sheaf and of a locally free sheaf so we may assume wrog that 
$W$ is locally free. Now using Serre duality, we may compute
$$
H^0(S,F^{n,*}_S(V^\vee)\otimes W)=H^1(S,F^{n,*}_S(V)\otimes W^\vee\otimes\Omega_{S/k_0})^\vee
$$
and the vector space $H^1(S,F^{n,*}_S(V)\otimes W^\vee\otimes\Omega_{S/k_0})$ vanishes for 
$n\gg 0$ because $V$ is cohomologically $p$-ample. \endProof

\begin{cor}
Suppose that 
$V$ is an ample bundle on $S$ and that $S$ is normal. Let $n_0\in\mN$ be such that 
$H^0(S,F^{n,*}_S(V^\vee)\otimes\Omega_{F_S})_g=0$ for all $n>n_0$. Let 
$S'$ be an irreducible scheme and let $\phi:S'\to S$ be a finite surjective morphism, which is generically inseparable. Then the map
$$
H^1(S,F^{n_0,*}_S(V^\vee))\to H^1(S',\phi^*(F^{n_0,*}_S(V^\vee)))
$$
is injective.
\label{cortor}
\end{cor}
\beginProof\, We may suppose that $S'$ is a normal scheme, since we may replace $S'$ by its normalization without restriction of generality. 
Let $H$ be the function field of $S$ and let $H'$ be the function field of $S'$. We are provided 
with a finite purely inseparable field extension $H\hookrightarrow H'$ and thus
for some sufficiently large $\ell_0\geqslant 1$, there exists a field extension $H'\hookrightarrow H$ and a commutative diagram

\vskip0.5cm
\centerline{
\begin{diagram}
H& \rInto& H'\\
\uTo^{=} & & \dInto\\
H&\rInto^{\,\,\,x^{p^{\ell_0}}}&H
\end{diagram}
}
\vskip0.5cm

Now notice that by Zariski's main theorem, $S'$ is canonically isomorphic to the normalisation of 
$S$ in $H'$. Similarly, the morphism $F_S^{\ell_0}:S\to S$ gives a presentation of $S$ as its own 
normalization in  $H$ via the $p^{\ell_0}$-th power map $H\stackrel{\,\,\,x^{p^{\ell_0}}}{\hookrightarrow} H$. There is thus a natural factorization $S\to S'\stackrel{\phi}{\to} S$, where the 
morphism $S\to S$ is given by  $F_S^{\ell_0}$. 
Using Corollary \ref{corinj}, we see that there is a natural 
injection $$H^1(S,F^{n_0,*}_S(V^\vee))\hookrightarrow H^1(S,F_S^{\ell_0,*}(F^{n_0,*}_S(V^\vee)))$$ and thus an injection $$
H^1(S,F^{n_0,*}_S(V^\vee))\to H^1(S',\phi^*(F^{n_0,*}_S(V^\vee))).
$$\endProof

\subsection{Jet schemes}

\label{ssec3}

The jet scheme construction described in \cite[sec. 2]{Rossler-MMML} provides a covariant functor
$$
\CY\mapsto J^i(\CY/U)
$$ 
from the category of quasi-projective schemes $\CY$ over $U$ to the category of quasi-projective 
schemes over $U$. 

The construction also provides an infinite tower of $U$-morphisms
$$
\dots\to J^2(\CY/U)\stackrel{\Lambda_{2,\CY}}{\to} J^1(\CY/U)\stackrel{\Lambda_{1,\CY}}{\to} \CY.
$$
If $Y$ is smooth over $U$, then the scheme $J^i(\CY/U)$, viewed as a $J^{i-1}(\CY/U)$-scheme via 
$\Lambda_{i,\CY}$, is a torsor under the 
vector bundle $\T\CY\otimes\Sym^i(\Omega_{U/k})$,  where $\T\CY:=(\Omega^{1}_{\CY/U})^\vee$. 
Here the bundle $\T\CY$ (resp.  $\Sym^i(\Omega_{U/k})$) is implicitly pulled back from $\CY$ (resp. $U$) to the scheme 
$J^{i-1}(\CY/U)$. 
 
 The functor $J^i(\cdot/U)$ preserves closed immersions and smooth morphisms. 
In particular, for any $i\in\mN$, there is a natural map 
\begin{equation}
\psi_{i,\CY}:\CY(U)=\Mor_U(U,\CY)\to 
\Mor_U(J^i(U/U),J^i(\CY/U))=\Mor_U(U,J^i(\CY/U)),
\label{fundlift}
\end{equation} and these maps are compatible with the 
morphisms $\Lambda_{i,\CY}$. 

Lastly, if a quasi-projective scheme $\CY$ over $U$ has the property that 
\mbox{$\CY\simeq \CY_0\times_k U$,} where $\CY_0$ is a quasi-projective scheme over 
$k$, then we are provided with a commutative diagram
\vskip0.5cm
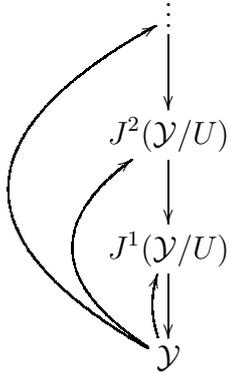
\centerline{
\xymatrix{ \vdots\ar[d]\\
      J^2(\CY/U)\ar[d] \\
      J^1(\CY/U)\ar[d]\\
\CY\ar@/^5pc/[uuu]\ar@/^3pc/[uu]\ar@/^/[u]}
}
\vskip0.5cm
In particular, in this situation, the scheme $J^i(\CY/U)$, viewed as a $J^{i-1}(\CY/U)$-torsor under 
under the vector bundle $\T\CY\otimes\Sym^i(\Omega_{U/k})$, is a trivial torsor. 
The curved arrows are functorial in an obvious manner with respect to the scheme $\CY_0/k$ and the 
isomorphism \mbox{$\CY\simeq \CY_0\times_k U$}. 

We refer to \cite[sec. 2]{Rossler-MMML} for the construction of jet schemes and 
for the proofs of the various assertions about them made above. 

\section{Proof of Theorem \ref{Mtheor}}

\label{sec3}

{\it Outline of the proof.} 

Suppose for the purposes of this outline that $X/K$ extends to a (not necessarily proper) smooth and quasi-projective $U$-scheme $\pi:\CX\to U$. 
First notice the following fact. If the conclusion of Theorem \ref{Mtheor} holds, then after possibly 
replacing $K$ by one of its finite separable extensions, 
there is an infinite commutative diagram

\vskip0.5cm
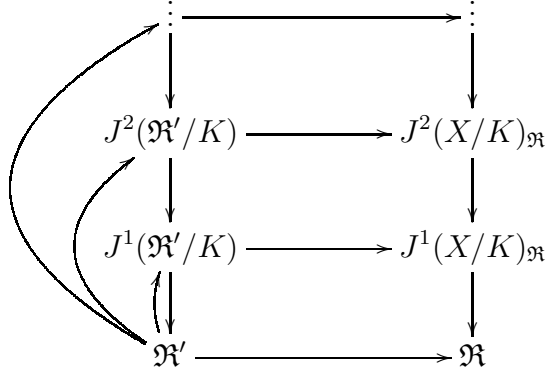
\begin{equation*}
\xymatrix{ \vdots\ar[rr]\ar[d]& & \vdots\ar[d]\\
      J^2({\mathfrak R}'/K)\ar[rr]\ar[d]  &   & J^2(X/K)_{\mathfrak R}\ar[d]\\
      J^1({\mathfrak R}'/K)\ar[rr]\ar[d]  &  & J^1(X/K)_{\mathfrak R}\ar[d]\\
{\mathfrak R}'\ar[rr]\ar@/^5pc/[uuu]\ar@/^3pc/[uu]\ar@/^/[u] & & {\mathfrak R}}
\end{equation*}
\vskip0.5cm

where ${\mathfrak R}'\simeq {\mathfrak R}'_0\times_k K$ and $J^i(X/K)_{\mathfrak R}=J^i(X/K)\times_X{\mathfrak R}. $
Here the horizontal arrows arise from the arrow ${\mathfrak R}'\to{\mathfrak R}$ by functoriality. 
The curved arrows arise from the isomorphism \mbox{${\mathfrak R}'\simeq {\mathfrak R}'_0\times_k K$} (see subsection \ref{ssec3}). 

In particular, we have an infinite commutative diagram
\vskip0.5cm
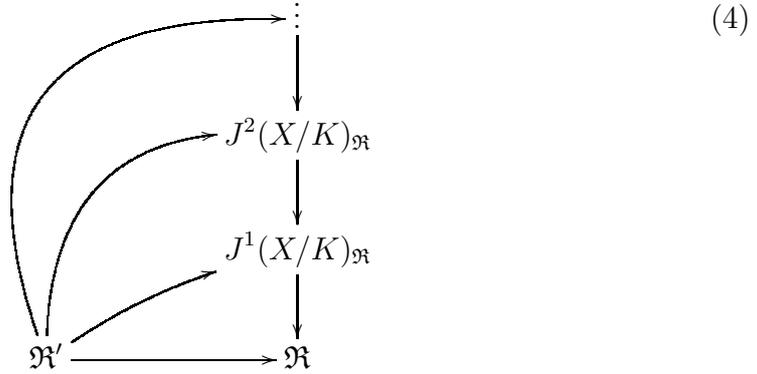
\begin{equation}
\xymatrix{ & & \vdots\ar[d]\\
        &   & J^2(X/K)_{\mathfrak R}\ar[d]\\
        &  & J^1(X/K)_{\mathfrak R}\ar[d]\\
{\mathfrak R}'\ar[rr]\ar@/^5pc/[uuurr]\ar@/^3pc/[uurr]\ar@/^/[urr] & & {\mathfrak R}}
\label{diagdiag}
\end{equation}
\vskip0.5cm

Step IV of the proof below will show that the existence of a diagram \refeq{diagdiag} is actually {\it equivalent} to 
the conclusion of Theorem \ref{Mtheor}. The method used in Step IV is analytic. 
Using the curved arrows in diagram \ref{diagdiag} and the definition of jet schemes we first construct a formal morphism \mbox{${\mathfrak R}'\to \CX$,} ie 
a morphism from ${\mathfrak R}'$ to $\CX$ viewed as formal schemes over the completion of 
$U$ along a suitable closed point. We then show using Grothendieck's formal GAGA theorem 
that this morphism has an algebraic model. 

In view of Step IV, it is thus sufficient to provide a diagram \refeq{diagdiag}. To do this, 
we use the fact that there are natural lifting maps $\psi_i:\CX(U)\to J^i(\CX/U)$ (see subsection \ref{ssec3}). 
Using a N\'eron desingularisation to construct $\CX$, we may assume that the natural map 
$\CX(U)\to X(K)$ is surjective so that $\CX(U)$ can be assumed to be large. 
The idea is now to consider the Zariski closure of the image of the map $\psi_i$ and to try to show that 
the base change to $\Spec(K)$ of one of these irreducible components qualifies as the image of a curved arrow ${\mathfrak R}'\to J^i(X/K)_{\mathfrak R}.$ 
We will see in Step II that there is a natural geometric candidate for the image of this arrow: the irreducible 
component $Z_i$ of positive dimension of the scheme $I({\phi_{n_i}})$ - see below.  
The crux of the proof is to show that the image of the lifting map $\CX(U)\to J^i(\CX/U)(U)$ meets 
$Z_i$ in a Zariski dense subset. To show this, we first compactify over $k$ all the objects 
in sight (ie the scheme $U$, the jet schemes $J^i(\CX/U)$ etc.), in order to 
avail ourselves of the Hilbert schemes of morphisms from the compactification of $U$ to 
the compactifications of the jet schemes. This is what we do in Step I, where these compactifications 
are constructed as canonically as possible by using the natural compactifications of torsors described in subsection \ref{ssec1}.  

In Step II, a height argument based on the existence of Hilbert schemes is then used to show that the schemes $I({\phi_{n_i}})$ 
must meet 'most' of the image of $\psi_i$ (see in particular Lemma \ref{hilblem} and Corollary \ref{corhilb}). The proof is made 
somewhat messy by the fact that some rational maps have to be made into 
morphisms in an arbitrary way (see especially after Lemma \ref{mmeq}). 

Once it is shown that the schemes $Z_i$ have the required properties, it remains 
to show that the morphisms between them eventually become birational (otherwise the 
diagram \refeq{diagdiag} cannot exist). This is what is shown in Step III. 
It is not difficult to show that the morphisms between the schemes $Z_i$ 
are generically purely inseparable but to show that they eventually become 
birational, we use the cohomological results of subsection \ref{ssec2}, especially Corollary \ref{corinj}. 
Here the torsor structure of the jet schemes is used and plugged into a cohomological machinery.

{\it The proof.}

We now suppose that the assumptions of Theorem \ref{Mtheor} are in force.  
We may suppose wrog that ${\mathfrak R}$ has positive dimension. Next, we may suppose wrog that $U$ is proper and smooth over $k$, since $U$ has dimension one. 

{\bf Step I.} {Compactifications.} 

By a compactification $\bar T$ of an $S$-scheme $T$, we shall mean a proper scheme $\bar T\to S$, which 
comes with an open immersion $T\hookrightarrow\bar T$ with dense image. 

Now choose any projective model $\bar\CX_0\to U$ of 
$X$ (this may for instance be obtained by embedding $X$ into some projective space over $K$). By applying N\'eron desingularization to $\bar\CX_0$ 
(see \cite[chap. 3, th. 2]{Bosch-Raynaud-Neron})), we obtain 
another projective model $\bar\CX_{00}$ of $X$ over $U$, with the property that the injection $\bar\CX_{00}^\sm(U)\hookrightarrow 
\bar\CX_{00}(U)=X(K)$ is a bijection. Here $\bar\CX_{00}^\sm\subset\bar\CX_{00}$ is the largest open subscheme $\bar\CX_{00}^\sm$of $\bar\CX_{00}$, such that $\bar\CX_{00}^\sm\to U$ is smooth. We now define $\CX:=\bar\CX^\sm_{00}$ and $\bar\CX:=\bar\CX_{00}$. 

To summarize: $\CX$ is a smooth model of 
$X$ over $U$ such that the natural map $\CX(U)\to X(K)$ is a bijection and we have an open immersion of $U$-schemes $\CX\hookrightarrow\bar\CX$, where 
$\bar\CX$ is projective over $U$. In particular the map $\CX\to U$ is surjective, since 
$\mathfrak R\not=\emptyset.$

We now choose specific compactifications  $\bar J^{i}(\CX/U)$ over $U$ for the jet schemes  $J^{i}(\CX/U)$. 

For $i=0$, we let  $\bar J^{0}(\CX/U)=\bar \CX$. We shall define the schemes $\bar J^{i}(\CX/U)$ inductively 
for $i\geqslant 0$. 

So suppose that the compactification  $\bar J^{i}(\CX/U)$ has already 
been constructed. As explained in subsection \ref{ssec3}, the $J^{i}(\CX/U)$-scheme 
$J^{i+1}(\CX/U)$  is a torsor under $F_i^\vee$, where \mbox{$F_i:=(\T \CX\otimes\Sym^{i+1}(\Omega_{U/k}))^\vee$} 
(viewed as a vector bundle over $J^{i}(\CX/U)$). Let 
\begin{equation}
0\to\CO_{J^{i}(\CX/U)}\to E_i\to F_i\to 0
\label{ext1}
\end{equation}
be an extension (unique up to non-unique isomorphism) associated with the class of 
$J^{i+1}(\CX/U)$ in $H^1(J^{i}(\CX/U),F_i^\vee)$ (see (xi) in subsection \ref{ssec1}). It was explained in (ii) subsection \ref{ssec1} that the 
$J^{i}(\CX/U)$-scheme $J^{i+1}(\CX/U)$ can be realized as the complement 
$\mP(E_i)\backslash\mP(F_i)$. 
We now define the 
compactification $\bar J^{i+1}(\CX/U)$ to be some  $\bar J^{i}(\CX/U)$-compactification of $\mP(E_i)$, 
such that the diagram 
\begin{diagram}
\mP(E_i) & \rInto & \bar J^{i+1}(\CX/U)\\
\dTo              &           &  \dTo\\
J^{i}(\CX/U)  & \rInto & \bar J^{i}(\CX/U)
\end{diagram}
is cartesian. This is possible because we may extend $E_i$ to a coherent sheaf $\bar E_i$ on 
$\bar J^{i}(\CX/U)$ and define $\bar J^{i+1}(\CX/U):=\mP(\bar E_i)$. 
We call $\bar\Lambda_{i+1}:\bar J^{i+1}(\CX/U)\to \bar J^{i}(\CX/U)$ the corresponding 
morphism.

The following diagram summarizes the resulting geometric configuration:
\begin{diagram}
J^{i+1}(\CX/U) & \rInto & \mP(E_i) & \rInto & \bar J^{i+1}(\CX/U)\\
\dTo^{\Lambda_{i+1}}              &           &  \dTo       &            &     \dTo_{\bar\Lambda_{i+1}}\\
J^{i}(\CX/U) & = &J^{i}(\CX/U)  & \rInto & \bar J^{i}(\CX/U)
\end{diagram}
Here the hooked horizontal arrows are open immersions and the square on the right is cartesian.

Recall the following key properties. The scheme $U$ is proper over $k$ and 
the schemes $\bar J^{i}(\CX/U)$ are proper over $U$. The morphisms 
$\bar J^{i+1}(\CX/U)\to \bar J^{i}(\CX/U)$ are proper. The schemes $J^{i}(\CX/U)$ and 
$\mP(E_i)$ are smooth and surjective onto $U$. 

Now remember that for each $i\geqslant 0$, there is a natural map \mbox{$\psi_i:\CX(U)\to J^{i}(\CX/U)(U)$,}
which is a section of the projection map $J^{i}(\CX/U)(U)\to\CX(U)$ (see 
subsection \ref{ssec3}). Abusing notation, 
we shall also call $\psi_i$ the map $\psi_i$ composed with the open immersion $J^{i}(\CX/U)\hookrightarrow \bar J^{i}(\CX/U).$

Finally, define $J^i(X/K):=J^i(\CX/U)_K$ and $\bar J^i(X/K):= \bar J^i(\CX/U)_K$.

{\bf Step II.} The schemes $Z_i\hookrightarrow J^i(X/K)$.

{\it We shall inductively construct closed integral subschemes $Z_i\hookrightarrow J^i(X/K)$ 
with the following properties. They are sent 
onto each other by the morphisms $\Lambda_{i,K}$. The morphisms $Z_{i+1}\to Z_{i}$ are finite, surjective  
and generically radicial and $Z_i$ is proper over $K$. In particular, $Z_i$ is closed in $\bar J^{i}(X/K)$. Finally, the image of $\CX(U)=X(K)$ by $\psi_{i,K}$ in 
$J^{i}(X/K)$ meets $Z_i(K)$ in a dense subset of $Z_i$ and $Z_0={\mathfrak R}$.}

The schemes $Z_i$ are defined via the following inductive procedure. 

We define $Z_0:={\mathfrak R}$. 

To define $Z_{i+1}$ from $Z_i$ notice that by Step I, we have an identification $\bar J^{i+1}(\CX/U)_{Z_{i}}=\mP(E_{i,Z_{i}})$. 

Notice also that $F_{i,Z_{i}}$ is ample (over $K$) since $Z_{i}\to X$ is finite. Thus by (vi) in subsection \ref{ssec1}, there exists 
a $K$-morphism 
$$
\phi_{n_i}:\mP(E_{i,Z_{i}})\to\mP^{n_i}_K,
$$
for some $n_i\in\mN$ that we fix. Using the notation of (viii) in subsection \ref{ssec1}, call $I({\phi_{n_i}})$ the union of the positive dimensional fibres of $\phi_{n_i}$. 

{\bf Substep II.1.} The closed set $I({\phi_{n_i}})$ dominates $Z_i$.

Write $H_i\subseteq \mP^{n_i}_K$ for a hyperplane such that $\phi_{n_i}^{-1}(H_i)=
\mP(F_{i,Z_{i}})$. To find such a hyperplane, consider that by (i) and (vi) in subsection \ref{ssec1}, the section of $\CO(F_{i,Z_{i}})^{\otimes(n_i+1)}$ corresponding 
to the Cartier divisor $(n_i+1)F_{i,Z_{i}}$ gives a section of $\CO_{\mP^{n_i}_K}(n_i+1)$. By construction, the hyperplane 
$H_i$ can be taken to be the reduced scheme underlying the zero scheme of this section.

Let $\bar\mP(E_{i,Z_{i}})$ be the Zariski closure 
of $\mP(E_{i,Z_{i}})$ in $\bar J^{i+1}(\CX/U)$ and let $\bar\mP(F_{i,Z_{i}})$ be the Zariski closure 
of $\mP(F_{i,Z_{i}})$ in $\bar J^{i+1}(\CX/U)$.

Now call $\Sigma_i\subseteq \CX(U)=X(K)$ the set of sections  $\sigma\in\CX(U)$ such that $\psi_{i,K}(\sigma)\in Z_i(K)$. 

\begin{lemma}
Let   
$\sigma\in\Sigma_i$ . We have 
$$
\psi_{i+1}(\sigma)\in \bar\mP(E_{i,Z_{i}})(U)
$$
and
\begin{equation}
\psi_{i+1}(\sigma)(U)\cap \bar\mP(F_{i,Z_{i}})=\emptyset.
\end{equation} 
\label{mmeq}
\end{lemma}
\beginProof\, The first equation follows from the definitions. 
The second equation follows from the fact that $\psi_{i+1}(\sigma)(U)\subseteq J^{i+1}(\CX/U)=\mP(E_i)\backslash\mP(F_i)$ and from the fact that we have 
a set-theoretic identity $\bar\mP(F_{i,Z_{i}})\cap 
\mP(E_i)\subseteq\mP(F_{i})$, since $\mP(F_{i})$ is a closed subset of $\mP(E_i)$. 
\endProof 

Now choose a proper birational U-morphism $b_i:\bar\mP'(E_{i,Z_{i}})\to \bar\mP(E_{i,Z_{i}})$, which is an isomorphism over $K$ and such 
that there exists a proper $U$-morphism $$\wt{\phi}_i:\bar\mP'(E_{i,Z_{i}})\to \mP^{n_i}_U,$$ with the property 
that $\phi_{n_i}\circ b_{i,K}=\wt{\phi}_{i,K}$. 

Let $\bar\mP'(F_{i,Z_{i}})$ be reduced Zariski closure of $b_{i,K}^{-1}(\mP(F_{i,Z_{i}}))$ in $\bar\mP'(E_{i,Z_{i}}).$ 
Let $\bar H_i\subseteq \mP^{n_i}_U$ be the Zariski closure of $H_i$. 

Now consider the Stein factorisation (see \cite[III, Cor. 11.5]{Hartshorne-Algebraic} for this)
$$
\bar\mP'(E_{i,Z_{i}})\stackrel{\rho_i}{\to} \bar\mP''(E_{i,Z_{i}})\stackrel{\bar\nu_i}{\to} \mP^{n_i}_U
$$
of $\wt{\phi}_i$. Notice that by construction the morphism $\rho_i$ is birational and the morphism $\bar\nu_i$ is finite.

The following commutative diagram sums up the geometric situation:

\vskip0.5cm
\centerline{
\xymatrix{
   &             \bar\mP(F_{i,Z_i})\ar@{^{(}->}[dr]                 & \bar J^{i+1}(\CX/U) & \bar J^{i+1}(X/K)\ar[l] & \mP(E_i)_K\ar[l]\\
 \bar\mP'(F_{i,Z_{i}})\ar@{^{(}->}[r]\ar[ur]\ar@/_4pc/[ddrr]& \bar\mP'(E_{i,Z_{i}})\ar[d]^{\rho_i}\ar[r]^{b_i}\ar[dr]^{\wt{\phi}_i}& \bar\mP(E_{i,Z_{i}})\ar@{^{(}->}[u] & \mP(E_{i,Z_{i}}) \ar[l]\ar[d]^{\phi_{n_i}}\ar@{^{(}->}[ur]&\mP(F_{i,Z_i})=\phi_{n_i}^{-1}(H_i)\ar[ddl]\ar[l] \\
 &             \bar\mP''(E_{i,Z_{i}})\ar[r]^{\bar\nu_i}                             & \mP^{n_i}_U & \mP^{n_i}_K\ar[l]& \\
   &                                        & \bar H_i \ar[u]        & H_i\ar[u]\ar[l]     &
                                           }
}
\vskip0.5cm

By the valuative criterion of properness, there is a natural map  $$\psi'_{i+1}:\CX(U)\to \bar\mP'(E_{i,Z_{i}})(U),$$ such that 
$b_i\circ\psi'_{i+1}=\psi_{i+1}$. 
\begin{lemma} 
There exists a constant $\beta_{i+1}\geqslant 0$, which is independent of $\sigma$, such that for all $\sigma\in\Sigma_i$, we have 
\begin{equation}
\length(\rho_i(\psi'_{i+1}(\sigma))(U)\cap\bar \nu_i^*(\bar H_i))\leqslant \beta_{i+1}.
\label{mmmeq}
\end{equation}
\end{lemma} 
Here $\cap$ refers to the scheme-theoretic intersection.
\beginProof\ 
Notice that $\bar\nu_i^{*}(\bar H_i)$ has a finite number of irreducible components. Among those, the only horizontal (ie 
dominating $U$) irreducible component 
is $\rho_i(\bar\mP'(F_{i,Z_{i}}))$. Furthermore, we have $\psi'_{i+1}(\sigma)\cap\bar\mP'(F_{i,Z_{i}})=\emptyset$ by equation \refeq{mmeq}. 
Thus 
$\rho_i(\psi'_{i+1}(\sigma))$ meets only the vertical irreducible components of $\bar\nu_i^{*}(\bar H_i)$. To conclude, notice that 
the intersection multiplicity of 
$\rho_i(\psi'_{i+1}(\sigma))$ with a fixed vertical component of $\bar\nu_i^{*}(\bar H_i)$ can be bounded independently of $\sigma$. 
\endProof

{\bf N.B.} The fact that we use N\'eron desingularisations is what makes it possible to prove the existence of 
the bound $\beta_{i+1}$ given in Lemma \ref{mmmeq}. Indeed, the fact that we use a N\'eron desingularisation 
implies that the Zariski closure in $\bar J^{i+1}(\CX/U)$ of the lifting of an element of 
$X(K)$ to $J^{i+1}(X/K)$ will lie over $J^{i}(\CX/U)$ and the restriction 
$\bar J^{i+1}(\CX/U)|_{J^{i}(\CX/U)}$ is the natural compactification $\mP(E_i)$, 
so that \refeq{mmeq} can be made to hold. 

Now we have 
\begin{lemma}
There exist a scheme $M_i$, which is quasi-projective over $k$ and a $U$-morphism 
$\mu_i:M_i\times_k U\to \bar\mP''(E_{i,Z_{i}})$, with the following properties:
\begin{itemize}
\item[$\bullet$] 
for all $P\in M_i(k)$, we have
$$
\deg((\mu_i\circ(P\times\Id_U))^*(\bar\nu_i^*(\CO(\bar H_i)))\leqslant \beta_{i+1};
$$
\item[$\bullet$] for any $U$-morphism $\kappa:U\to \bar\mP''(E_{i,Z_{i}})$ such that 
$$
\deg((\mu_i\circ (P\times\Id_U))^*(\bar\nu_i^*(\CO(\bar H_i)))\leqslant \beta_{i+1}
$$
there is a $P\in M_i(k)$ such that $\kappa=\mu_i\circ(P\times\Id_U)$. 
\end{itemize}
Here $\beta_{i+1}$ is the constant appearing in Lemma \ref{mmmeq}.
\label{hilblem}
\end{lemma}
\beginProof\  
Notice that the morphism $\bar\nu_i$ is finite and thus the divisor 
$\bar\nu_i^*(\bar H_i)$ is ample. The existence of $M_i$ is now a consequence of the theory of 
Hilbert schemes.
\endProof

\begin{cor}
For almost all the sections $\sigma\in\Sigma_i$, we have 
$$\psi_{i+1,K}(\sigma)\in (I({\phi_{n_i}}))(K).$$
\label{corhilb}
\end{cor}
\beginProof\,
Let $M_i$ be as in Lemma \ref{hilblem}. We may suppose wrog that $\Sigma_i$ is infinite. 
In view of \refeq{mmmeq}, there is a natural map $\gr:\Sigma_i\to M_i(k)$, such 
that $$(\mu_i\circ(\gr(\sigma)\times\Id_U))=\sigma$$ for all $\sigma\in\Sigma_i$.

On the other hand, consider 
$N$, a reduced subscheme of $M_i.$ We have by construction a $K$-morphism 
$N_K\to \bar\mP''(E_{i,Z_{i}})_K$. Let 
$U'''_i\subseteq\mP^{n_i}_K$ be the locus of points where the fiber of 
$\wt{\phi}_{i,K}$ is finite and let $U''_i:=\bar\nu_i^{-1}(U'''_i)$. Note that by the construction of the 
Stein factorisation, the morphism $\rho_i|_{U''_i}:\rho_i^{-1}(U''_i)\to U''_i$ is an isomorphism. 
If the conclusion of the corollary fails, we may thus find $N$ as above such that  $N$ is integral and $\dim(N)>0$, such that 
$N_K$ intersects $U''_i$ and such that $N\cap\gr(\Sigma_i)$ is dense in $N$. 

We suppose that such an $N$ exists to obtain a contradiction. Restricting the size of $N$, we may then suppose that 
the image of $N_K\to \bar\mP''(E_{i,Z_{i}})_K$ lies inside $U''_i$ 
and thus we obtain a morphism $l_N:N_K\to\mP(E_{i,Z_{i}})\subseteq \bar J^{i+1}(X/K)$.

Consider that we also have by functoriality a morphism 
$$J^{i+1}(N_K/K)\to J^{i+1}(X/K)$$
 arising from the morphism $N_K\to X$, which is the  
composition of $l_N$ with the projection $\mP(E_{i,Z_i})\to X$. 
The morphism $J^{i+1}(N_K/K)\to J^{i+1}(X/K)$ can be composed with the natural section 
$N_K\to J^{i+1}(N_K/K)$ (see subsection \ref{ssec3}) 
to obtain a second $K$-morphism $$l'_N:N_K\to J^{i+1}(X/K)\subseteq \bar J^{i+1}(X/K).$$ By construction, the morphisms $l_N$ and $l'_N$ coincide 
on a dense set of $K$-points and thus $l_N=l'_N$. 
 
Now let $C\subseteq N$ be a smooth curve. Repeating the construction of $l'_N$ for $C$, 
we obtain a morphism $l'_C:C_K\to J^{i+1}(X/K)$ and furthermore the composition of $l'_N=l_N$ with $C\to N$ 
is  $l'_C$. Now consider a smooth compactification $\bar C$ of $C$ over $k$. By the valuative 
criterion of properness, there is a unique $K$-morphism $\bar C_K\to X$ extending 
the morphism $C_K\to X$. Following the steps 
of the construction of $l'_C$, we again obtain 
a $K$-morphism $$l'_{\bar C}:\bar C_K\to J^{i+1}(X/K),$$ which extends $l'_C$. On the other hand, 
$l'_{\bar C}(\bar C_K)\subseteq \mP(E_{i,Z_i})\backslash\mP(F_{i,Z_i})$ by construction and 
the image 
$\phi_{n_i}(l'_{\bar C}(\bar C_K))$ is closed in $\mP^{n_i}_K$. Furthermore, by the definition 
of $H_i$, $$\phi_{n_i}(l'_{\bar C}(\bar C_K))\cap H_i=\emptyset.$$ Thus $\phi_{n_i}(l'_{\bar C}(\bar C_K))$ is the underlying set of an integral scheme, which is affine and proper over $K$ and 
thus consists of 
a closed point. Since by construction $l'_C(C_K)\not\subseteq I({\phi_{n_i}})$, we see that 
$l'_C(C_K)$ is a closed point. Now any two closed points in the smooth locus $N^\sm$ of $N$ can be connected 
by a smooth curve lying in $N^\sm$ (see for instance \cite{Katz-Lefschetz}). Hence the image of $N^\sm(k)$ by $l'_N$ is a closed point $P$.
Since the image of $N^\sm(k)$ in $N^\sm(K)$ is dense $N_K$ the fibre $(l'_N)^{-1}(P)$ is dense in $N_K$ and thus 
$l'_N(N_K)=P$. This contradicts the fact that $N\cap\gr(\Sigma_i)$ is dense in $N$ (because 
the elements of $\Sigma_i$ correspond to pairwise distinct points in $X(K)$). 
 \endProof

We can now conclude that $I({\phi_{n_i}})$ dominates $Z_{i}$, since by Corollary \ref{corhilb} the image of $I({\phi_{n_i}})$ in $Z_i$ contains the image 
by $\psi_{i,K}$ of a cofinite set of $\Sigma_i$ and this set is dense in $Z_i$ by assumption.

{\bf Substep II.2.} Definition of $Z_{i+1}$.

Since $\psi_{i,K}(\Sigma_i)$ meets $Z_{i}$ in a dense subset by assumption, we can now deduce from (ix) in subsection \ref{ssec1} that 
$I({\phi_{n_i}})$ has a single irreducible component of positive dimension, which is finite, surjective and generically radicial over $Z_{i}$. 

{\it We now define $Z_{i+1}$ to be this irreducible component.}

We have shown that $Z_{i+1}$ has all the required properties.

{\bf Step III.} {Purely inseparable trivialization}.

Let $\wt{Z}_0\to Z_0$ be the normalisation of $Z_0$. We note the important fact that $\wt{Z}_0$ is geometrically 
reduced over $K$ because $\wt{Z}_0(K)$ is dense in $\wt{Z}_0$ (see eg \cite[Lemma 1]{AV-Toward} for this). For every $i\geqslant 0$, write $\wt{Z}_{i}\to\wt{Z}_0$
 for the $\wt{Z}_0$-scheme obtained from $Z_i$ by base-change. Similarly, write $\wt{J}^i(X/K)$ for the pull-back of $J^i(X/K)$ to $\wt{Z}_0$. We denote by $\wt{F}_0$ the pull-back to $\wt{Z}_0$ of the vector bundle $F_0\simeq(\T X\otimes\Omega_{K/k})^\vee.$

{\it Suppose that $\charac(k)=0$.} 

The scheme $\wt{J}^1(X/K)$ is by construction a torsor under $\wt{F}^\vee_0.$
Furthermore, the morphism $\wt{Z}_1^\red\to \wt{Z}_0$ is an isomorphism by 
Zariski's main theorem and it trivializes the torsor $\wt{J}^1(X/K)$ by construction.

Repeating this reasoning for $\wt{Z}_2^\red$ over $\wt{Z}_1^\red\simeq \wt{Z}_0$, $\wt{Z}_3^\red$ over 
$\wt{Z}_2^\red$ and so forth, 
we see that the natural morphisms $\wt{Z}_{i+1}^\red\to \wt{Z}_{i}^\red$ are isomorphisms for all $i\geqslant 0$.

{\it Now suppose until the end of step III that $\charac(k)>0$.}

Write $F^{m_0,*}_{\wt{Z}_0} \wt{Z}_i$ for the base-change of $\wt{Z}_i\to \wt{Z}_0$ by $F^{m_0}_{\wt{Z}_0}$. Let $$\pi_{i,m_0}:(F^{m_0,*}_{\wt{Z}_0} \wt{Z}_i)^\red\to \wt{Z}_0$$ be the natural morphism. 

Similarly, write 
$F^{m_0,*}_{\wt{Z}_0}\wt{J}^i(X/K)$ for the base-change of $\wt{J}^i(X/K)$ by $F^{m_0}_{\wt{Z}_0}$. 

We now define
$$
m_0:=\sup\{m\in\mN|H^0(\wt{Z}_0,F_{\wt{Z}_0}^{m,*}(\wt{F}_0^\vee)\otimes\Omega_{F_{\wt{Z}_0}})_g\not=0\}.$$ By Proposition \ref{propn}, the quantity $m_0$ is finite.

The scheme $F^{m_0,*}_{\wt{Z}_0}\wt{J}^1(X/K)$ is by construction a torsor under $F^{m_0,*}_{\wt{Z}_0}(\wt{F}_0).$
Furthermore, the morphism $\pi_{1,m_0}$ trivializes this torsor by construction, since there is a $\wt{Z}_0$-morphism $(F^{m_0,*}_{\wt{Z}_0}\wt{Z}_1)^\red\to F^{m_0,*}_{\wt{Z}_0}\wt{J}^1(X/K)$. By the assumption on $m_0$ and because $\pi_{1,m_0}$ is finite and generically radicial, the $F^{m_0,*}_{\wt{Z}_0}(\wt{F}_0)$-torsor $F^{m_0,*}_{\wt{Z}_0}\wt{J}^1(X/K)$ is actually trivial (use Corollary 
\ref{cortor}). 
Let 
$$t:\wt{Z}_0\to F^{m_0,*}_{\wt{Z}_0}\wt{J}^1(X/K)$$ be a section. By (x) in subsection \ref{ssec1} we have that $t_*(\wt{Z}_0)=(F_{\wt{Z}_0}^{m_0,*}\wt{Z}_1)^\red$ and thus 
the morphism \mbox{$(F_{\wt{Z}_0}^{m_0,*}\wt{Z}_1)^\red\to \wt{Z}_0$} is an isomorphism. 

Repeating this reasoning for $(F^{m_0,*}_{\wt{Z}_0} \wt{Z}_2)^\red$ over $(F_{\wt{Z}_0}^{m_0,*}\wt{Z}_1)^\red\simeq \wt{Z}_0$, $(F^{m_0,*}_{\wt{Z}_0}\wt{Z}_3)^\red$ over $(F^{m_0,*}_{\wt{Z}_0}\wt{Z}_2)^\red$ and so forth, 
we see that the natural morphisms $$(F^{m_0,*}_{\wt{Z}_0} \wt{Z}_{i+1})^\red\to (F^{m_0,*}_{\wt{Z}_0}\wt{Z}_{i})^\red$$ are isomorphisms from all $i\geqslant 0$.

{\bf Step IV.} {Formalization and utilization of Grothendieck's GAGA.} 

In this fourth step, the argument will be written up under the assumption that $\charac(k)>0$. 
The argument goes through verbatim when $\charac(k)=0$, if one sets $m_0=0$ in the text below. 

We now choose a non empty affine open subset $U_0\subseteq U$ such that $\CX_{U_0}\to U_0$ is 
smooth and projective and we let $\CZ_0$ be the reduced Zariski closure of $Z_0={\mathfrak R}$ in $\CX_{U_0}$. 
We let $\wt{\CZ}_0$ be the normalisation of $\CZ_0$. By construction, the natural morphism $\wt{\CZ}_0\to\CZ_0$ is finite 
and birational and the morphism $\wt{\CZ}_{0,K}\to \CZ_{0,K}=Z_0$ can be identified 
with the morphism $\wt{Z}_0\to Z_0$. 

We may assume wrog that the natural map
\begin{equation}
H^1(\wt{\CZ}_0,F^{m_0,*}_{\wt{\CZ}_0}(\T \CX_{U_0}\otimes\Omega_{U_0/k}))\to H^1(\wt{Z}_0,F^{m_0,*}_{\wt{Z}_0}(\T X\otimes\Omega_{K/k}))
\label{RMF}
\end{equation}
is injective (here as before, $\T \CX_{U_0}$ and $\Omega_{U_0/k}$ are identified with their pull-backs 
to $\wt{\CZ}_0$ and $\T X$ and $\Omega_{K/k}$ are identified with their pull-backs to $\wt{Z}_0$). 

Indeed, let $\pi_{\wt{\CZ}_0\to U_0}:\wt{\CZ}_0\to U_0$ and $\pi_{\wt{Z}_0\to\Spec(k)}:\wt{Z}_0\to\Spec(k)$ be the 
structural morphisms. Since $U_0$ is affine, we have 
$$
H^1(\wt{\CZ}_0,F^{m_0,*}_{\wt{\CZ}_0}(\T \CX_{U_0}\otimes\Omega_{U_0/k}))=
H^0(U_0,{\rm R}^1\pi_{\wt{\CZ}_0\to U_0,\,\ast}(F^{m_0,*}_{\wt{\CZ}_0}(\T \CX_{U_0}\otimes\Omega_{U_0/k})))
$$
and the map $H^1(\wt{\CZ}_0,F^{m_0,*}_{\wt{\CZ}_0}(\T \CX_{U_0}\otimes\Omega_{U_0/k}))\to H^1(\wt{Z}_0,F^{m_0,*}_{\wt{Z}_0}(\T X\otimes\Omega_{K/k}))$ is induced by the morphism of sheaves
$$
{\rm R}^1\pi_{\wt{\CZ}_0\to U_0,\,\ast}(F^{m_0,*}_{\wt{\CZ}_0}(\T \CX_{U_0}\otimes\Omega_{U_0/k}))
\to {\rm R}^1\pi_{\wt{Z}_0\to\Spec(k),\,\ast}(F^{m_0,*}_{\wt{Z}_0}(\T X\otimes\Omega_{K/k}))
$$
which is injective if the sheaf ${\rm R}^1\pi_{\wt{\CZ}_0\to U_0,\,\ast}(F^{m_0,*}_{\wt{\CZ}_0}(\T \CX_{U_0}\otimes\Omega_{U_0/k}))$ is torsion free. This 
will be the case if $U_0$ is sufficiently small, since the sheaf \mbox{${\rm R}^1\pi_{\wt{\CZ}_0\to U_0,\,\ast}(F^{m_0,*}_{\wt{\CZ}_0}(\T \CX_{U_0}\otimes\Omega_{U_0/k}))$} is a coherent sheaf.

We write $\wt{J}^i(\CX_{U_0}/U_0)$ for the base change of $J^i(\CX_{U_0}/U_0)$ to $\wt{\CZ}_0$ and 
$F^{m_0,*}_{\wt{\CZ}_0}\wt{J}^i(\CX_{U_0}/U_0)$ for the base change of $\wt{J}^i(\CX_{U_0}/U_0)$ by $F^{m_0}_{\wt{\CZ}_0}$.

Because of the existence of the injection \refeq{RMF} and the end of Step III, the \mbox{$F^{m_0,*}_{\wt{\CZ}_0}(\T \CX_{U_0}\otimes\Omega_{U_0/k})$} torsor
$F^{m_0,*}_{\wt{\CZ}_0}\wt{J}^1(\CX_{U_0}/U_0)$ is trivial. Furthermore, if we choose a section $$\wt{\CZ}_0\hookrightarrow F^{m_0,*}_{\wt{\CZ}_0}\wt{J}^1(\CX_{U_0}/U_0)$$ then by (x) in subsection \ref{ssec1} and again by \refeq{RMF}, the pull-back by this section of the 
\mbox{$F^{m_0,*}_{\wt{\CZ}_0}(\T \CX_{U_0}\otimes\Omega_{U_0/k})$} torsor 
$F^{m_0,*}_{\wt{\CZ}_0}\wt{J}^2(\CX_{U_0}/U_0)$ is also trivial. Continuing in this way, we get 
a sequence of sections $\wt{\CZ}_0\hookrightarrow F^{m_0,*}_{\wt{\CZ}_0}\wt{J}^i(\CX_{U_0}/U_0)$. The following 
commutative diagram summarises the situation:   

\vskip0.5cm
\begin{equation*}
\xymatrix{ \vdots\ar[d]& & \vdots\ar[d]^{\sim}\ar[ll]\ar[r]& \vdots\ar[d]\\
      F^{m_0,*}_{\wt{\CZ}_0}\wt{J}^2(\CX_{U_0}/U_0)\ar[d]  &   & 
      (F^{m_0,*}_{\wt{Z}_0} \wt{Z}_{2})^\red\ar[d]^{\sim}\ar[r]\ar[ll] & \wt{Z}_2\ar[d]\\
      F^{m_0,*}_{\wt{\CZ}_0}\wt{J}^1(\CX_{U_0}/U_0)\ar[d]  &  & (F^{m_0,*}_{\wt{Z}_0} \wt{Z}_{1})^\red\ar[d]^{\sim}\ar[ll]\ar[r]& \wt{Z}_1\ar[d]\\
\wt{\CZ}_0\ar@/^10pc/[uuu]\ar@/^6pc/[uu]\ar@/^/[u] & & \wt{Z}_0\ar[ll]\ar[r]^{F^{m_0}_{\wt{Z}_0}}& \wt{Z}_0
}
\end{equation*}
\vskip0.5cm

Now write $\wt{J}^i(\CZ_0/U_0)$ for the base-change of $J^i(\CZ_0/U_0)$ to $\wt{\CZ}_0$ and 
$F^{m_0,*}_{\wt{\CZ}_0}\wt{J}^i(\CZ_0/U_0)$ for the base-change of $\wt{J}^i(\CZ_0/U_0)$ by $F^{m_0}_{\wt{\CZ}_0}$. 
Similarly, write $\wt{J}^i(Z_0/K)$ for the base-change of $J^i(Z_0/K)$ to $\wt{Z}_0$ and 
$F^{m_0,*}_{\wt{Z}_0}\wt{J}^i(Z_0/K)$ for the base-change of $\wt{J}^i(Z_0/K)$ by $F^{m_0}_{\wt{Z}_0}$. 

Notice that we have a natural identification $\wt{J}^i(Z_0/K)\simeq \wt{J}^i(\CZ_0/U_0)_K$ and 
$$(F^{m_0,*}_{\wt{\CZ}_0}\wt{J}^i(\CZ_0/U_0))_K\simeq F^{m_0,*}_{\wt{Z}_0}\wt{J}^i(Z_0/K).$$

We claim that each of the sections $\wt{\CZ}_0\hookrightarrow F^{m_0,*}_{\wt{\CZ}_0}\wt{J}^i(\CX_{U_0}/U_0)$ factors through $F^{m_0,*}_{\wt{\CZ}_0}\wt{J}^i(\CZ_0/U_0)$ viewed as a closed 
subscheme of $F^{m_0,*}_{\wt{\CZ}_0}\wt{J}^i(\CX_{U_0}/U_0).$ 
Since $\wt{\CZ}_0$ is integral and dominates $U_0$, to establish this, it is sufficient 
to show that the image of the corresponding section $\wt{Z}_0\to F^{m_0,*}_{\wt{Z}_0}\wt{J}^i(X/K)$ is contained 
in the closed subscheme $$F^{m_0,*}_{\wt{Z}_0}\wt{J}^i(Z_0/K)\hookrightarrow F^{m_0,*}_{\wt{Z}_0}\wt{J}^i(X/K).$$

To see this, first notice that we have a commutative diagram
\vskip0.5cm
\begin{equation}
\xymatrix{
\vdots\ar@{^{(}->}[r]\ar[d] & \vdots\ar[d]\ar@{^{(}->}[r] & \vdots\ar[d]\\
Z_2\ar@{^{(}->}[r]\ar[d]  & J^2(Z_0/K)\ar@{^{(}->}[r]\ar[d] & J^2(X/K)\ar[d]\\
Z_1\ar@{^{(}->}[r] & J^1(Z_0/K)\ar@{^{(}->}[r]\ar[d] & J^1(X/K)\ar[d]\\
 & Z_0\ar@{^{(}->}[r] & X
}
\label{vipIII}
\end{equation}
\vskip0.5cm

This does {\it not} follow from the general properties of jet schemes but follows from the fact 
that for all $i\geq 0$, the set $\psi_i(Z_0(K))\cap Z_i(K)$ is (by construction) dense  in $Z_i$ 
and from the fact that $
\psi_i(Z_0(K)))\subseteq J^i(Z_0/K)(K)$ by the functoriality of jet schemes. Applying base-change to the commutative diagram \refeq{vipIII}, we deduce that there is a commutative diagram

\vskip0.5cm
\begin{equation*}
\xymatrix{
\vdots\ar@{^{(}->}[r]\ar[d]^{\sim} & \vdots\ar[d]\ar@{^{(}->}[r] & \vdots\ar[d]\\
(F^{m_0,*}_{\wt{Z}_0}\wt{Z}_2)^\red\ar@{^{(}->}[r]\ar[d]^{\sim}  & F^{m_0,*}_{\wt{Z}_0}\wt{J}^2(Z_0/K)\ar@{^{(}->}[r]\ar[d] & F^{m_0,*}_{\wt{Z}_0}\wt{J}^2(X/K)\ar[d]\\
(F^{m_0,*}_{\wt{Z}_0}\wt{Z}_1)^\red\ar@{^{(}->}[r]\ar[dr]^{\sim} & F^{m_0,*}_{\wt{Z}_0}\wt{J}^1(Z_0/K)\ar@{^{(}->}[r]\ar[d] & F^{m_0,*}_{\wt{Z}_0}\wt{J}^1(X/K)\ar[d]\\
 & \wt{Z}_0\ar[r]^{=} & \wt{Z}_0
}
\end{equation*}
\vskip0.5cm

and this establishes our claim.

Now choose a closed point $u_0\in U_0.$ View $u_0$ as a closed subscheme of $U_0$.  For any $i\geqslant 0$, let $u_i$ be the $i$-th infinitesimal neighborhood of $u_0\simeq \Spec\, k$ in $U_0$ (so that there is no ambiguity of notation for $u_0$). Notice that $u_i$ has a natural structure of $k$-scheme. We may assume wrog that the morphism $\wt{\CZ}_{0,{u_0}}\to\CZ_{u_0}$ is birational. Recall that by the definition of jet schemes (see 
\cite[sec. 2]{Rossler-MMML}), the scheme $J^i({\CZ}_0/U_0)_{u_0}$ represents the functor 
on $k$-schemes
$$
T\mapsto\Mor_{u_i}(T\times_k{u_i},{\CZ}_{0,u_i}).
$$
Thus the infinite chain of compatible morphisms $\wt{\CZ}_0\hookrightarrow F^{m_0,*}_{{\CZ}_0}\wt{J}^i({\CZ}_0/U_0)$ gives rise to  $u_i$ morphisms
\begin{equation}
\wt{\CZ}_{0,{u_0}}^{(p^{-m_0})}\times_{k_0} u_i\to {\CZ}_{0,u_i}
\label{formmor}
\end{equation}
compatible with each other under base-change. Here $\wt{\CZ}_{0,{u_0}}^{(p^{-m_0})}$ denotes 
the base-change of the $k$-scheme $\wt{\CZ}_{0,{u_0}}$ by the $m_0$-th power 
of the inverse of the Frobenius automorphism of $k$. 

Let $\widehat{U}_{u_0}$ be the completion of the local ring of $U$ at $u_0$. 
View the $\widehat{U}_{u_0}$-schemes $\wt{\CZ}_{0,u_0}^{(p^{-m_0})}\times_k \widehat{U}_{u_0}$ and ${\CZ}_{0,\widehat{U}_{u_0}}$ as formal 
schemes over $\widehat{U}_{u_0}$ in the next sentence. The family
of morphisms \refeq{formmor} provides us with a morphism of formal schemes 
$$
\wt{\CZ}_{0,{u_0}}^{(p^{-m_0})}\times_k \widehat{U}_{u_0}\to{\CZ}_{0,\widehat{U}_{u_0}}
$$
and since both schemes are projective over $\widehat{U}_{u_0}$, Grothendieck's formal GAGA theorem shows that this morphism of 
formal schemes comes from  a unique morphism of 
$\widehat{U}_{u_0}$-schemes
$$
\iota:\wt{\CZ}_{0,{u_0}}^{(p^{-m_0})}\times_k \widehat{U}_{u_0}\to{\CZ}_{0,\widehat{U}_{u_0}}.
$$
By construction, at the closed point $u_0$ of $\widehat{U}_{u_0}$, the morphism $\iota$ specializes to the morphism $F^{m_0}_{\wt{\CZ}_{0,u_0}}$ composed 
with the finite and birational morphism $\wt{\CZ}_{0,u_0}\to\CZ_{u_0}$. 
Now the set of points of $\CX_{\widehat{U}_{u_0}}$, where the fibres of $\iota$ are non-empty and of dimension $>0$ is closed (see \cite[IV,13.1.5]{EGA}). On the other hand, since $\CX_{\widehat{U}_{u_0}}$ is proper over $\Spec(\widehat{U}_{u_0})$ and since 
 this set does not meet the special fibre $\CX_{{u_0}}$, it must be empty. Thus  
the morphism $\iota$ is finite over the generic point of $\widehat{U}_{u_0}$. 

For use in Step V, we also 
record the fact that if $\iota$ can be shown to be flat then its degree $\deg(\iota)$ is equal to the degree $\deg(\iota_{u_0})$ of its base change to the special fibre  $\CX_{{u_0}}$, which 
is $p^{\dim({\mathfrak R})m_0}.$ ($\ast\ast$)

Let $\widehat{K}$ be the function field of $\widehat{U}_{u_0}$. Since $k$ is an excellent field, we know that the field extension $\widehat{K}|K$ is separable (see eg \cite[§8.2, Ex. 2.34]{Liu-Algebraic}). On the other hand 
the just constructed finite morphism $\wt{\CZ}_{0,{u_0}}^{(p^{-m_0})}\times_k\widehat{K}\to Z_{0,\widehat{K}}$ is defined 
over a finitely generated subfield $K'$ (as a field over $K$) of $\widehat{K}$. The field extension $K'|K$ is then still separable, 
so that by the theorem on separating transcendence bases, there exists a variety $U'/K$, which is smooth over $K$ and whose function field is $K'$. Furthermore, possibly replacing $U'$ by one of its open subschemes, we may assume that 
the morphism $\wt{\CZ}_{0,{u_0}}^{(p^{-m_0})}\times_k {K}'\to Z_{0,{K}'}$ extends to a finite  morphism 
$$\alpha:\wt{\CZ}_{0,{u_0}}^{(p^{-m_0})}\times_k U'\to Z_{0,U'}.$$  
Furthermore, we may assume wrog that $\alpha$ is flat if $\iota$ is flat. 

Let $P\in U'(K^\sep)$ be a $K^\sep$-point over $K$ (the set $U'(K^\sep)$ is not empty because 
$U'$ is smooth over $K$). The base change of the morphism $\alpha$ by $P$ gives a morphism 
$$\wt{\CZ}_{0,{u_0}}^{(p^{-m_0})}\times_k K^\sep\to Z_{0,K^\sep}.$$ This is the morphism $h$ advertised in Theorem \ref{Mtheor}. 
If $\iota$ is flat, then we may assume that $\deg(\iota)=\deg(h)$. 

{\bf Step V}. Degree estimates.

In this last step, we shall prove that if ${\mathfrak R}=X$ then 

- the morphism $\iota$ is flat;

- the inequality 
\begin{equation}
\charac(k)>\dim(X)^2\int_X{\rm c}_1(\Omega_X)^{\dim(X)}
\label{finEE}
\end{equation}
implies that $m_0=0$. 

According to the end of Step V (see in particular $(\ast\ast)$), if $\iota$ is flat 
and $m_0=0$ then 
the morphism $h$ can be assumed to be an isomorphism. This will thus conclude the proof 
of Theorem \ref{Mtheor}. 

So suppose until the end of Step IV that ${\mathfrak R}=X$.
 Since both the source and the target of the morphism 
$\iota$ are regular schemes, the morphism 
$\iota$ is automatically flat by 'miracle flatness' (see \cite[Th. 23.1]{Matsumura-Commutative}). Next recall that by definition 
$$
m_0:=\sup\{m\in\mN|H^0(\wt{Z}_0,F_{\wt{Z}_0}^{m,*}(\wt{F}_0^\vee)\otimes\Omega_{F_{\wt{Z}_0}})_g\not=0\}.
$$
Since $X={\mathfrak R}$ is smooth over $K$, we have $\wt{F}_0=\Omega_{X/K}$ and $F_{\wt{Z}_0}=F_X$. Furthermore, let $f:X\to\Spec(K)$ be the structural morphism. Notice that we have an exact sequence
$$
F_X^*(\Omega_{X/k})\stackrel{F_X^*}{\to} \Omega_{X^{(p)}/k}\to\Omega_{F_X}\to 0
$$
where $X^{(p)}$ is the $k$-scheme $X$ with the structural morphism $(F_{\Spec(k)})^{-1}\circ f$. 
For a function $x$ defined on an open set of $X$ we have $F_X^*(x)=x^p$, 
where $p:=\charac(K)$. Thus we may compute $$F_X^*(d(x))=d(F_X^*(x))=d(x^p)=p\cdot x^{p-1}\cdot d(x)=0$$ 
and thus the morphism $F_X^*:F_X^*(\Omega_{X/k})\stackrel{F_X^*}{\to} \Omega_{X^{(p)}/k}$ 
vanishes. So we have a natural isomorphism $\Omega_{X^{(p)}/k}\simeq\Omega_{F_X}.$
Note that we also an isomorphism $\Omega_{X^{(p)}/k}\simeq \Omega_{X/k}$, since 
 $(F_{\Spec(k)})^{-1}$ is an isomorphism.

We also 
have an exact sequence
$$
f^*(\Omega_{K/k})\to\Omega_{X/k}\to\Omega_{X/K}\to 0
$$
and since the morphism $f$ has a section (since $X(K)\not=\emptyset$ because $X(K)$ is dense in $X$!), this sequence is exact on the left and splits. Thus we have a non-canonical isomorphism
$$
f^*(\Omega_{K/k})\oplus\Omega_{X/K}\simeq\Omega_{F_X}
$$
In particular, $\Omega_{F_X}$ is locally free and thus we have
\begin{eqnarray*}
&&m_0=\sup\{m\in\mN|H^0(X,F_{X}^{m,*}(\Omega_{X/K}^\vee)\otimes\Omega_{F_{X}})\not=0\}=\\
&&\hskip-1.5cm\max\big(\,\sup\{m\in\mN|H^0(X,F_{X}^{m,*}(\Omega_{X/K}^\vee)\otimes f^*(\Omega_{K/k}))\not=0\},\, 
\sup\{m\in\mN|H^0(X,F_{X}^{m,*}(\Omega_{X/K}^\vee)\otimes\Omega_{X/K})\not=0\}\,\big)
\end{eqnarray*}
Now $f^*(\Omega_{K/k})\simeq\CO_X$ (non canonically) and for 
all $m\geq 0$, the sheaf 
$F_{X}^{m,*}(\Omega_{X/K})$ is ample. Thus the sheaf $F_{X}^{m,*}(\Omega_{X/K}^\vee)$ has no global sections. 
In particular, for all $m\geq 0$, we have $H^0(X,F_{X}^{m,*}(\Omega_{X/K}^\vee)\otimes f^*(\Omega_{K/k}))=0$ and thus we have 
\begin{eqnarray*}
m_0=\sup\{m\in\mN|H^0(X,F_{X}^{m,*}(\Omega_{X/K}^\vee)\otimes\Omega_{X/K})\not=0\}.
\end{eqnarray*}
Now it is a special case of \cite[Lemma-Def. 3.3]{Rossler-SSS} that we have
$$
(\charac(\bar K))^{\sup\{m\in\mN|H^0(X_{\bar K},F_{X_{\bar K}}^{m,*}(\Omega_{X_{\bar K}/\bar K}^\vee)\otimes\Omega_{X_{\bar K}/\bar K})\not=0\}}\leqslant{\bar\mu_\max(
\Omega_{X_{\bar K}/\bar K})\over 
\bar\mu_\min(\Omega_{X_{\bar K}/\bar K})}
$$
where $\bar\mu_\max(\Omega_{X_{\bar K}/\bar K})$ 
(resp. $\bar\mu_\min(\Omega_{X_{\bar K}/\bar K})$) is the Frobenius stabilised 
maximal (resp. minimal) slope of $\Omega_{X_{\bar K}/\bar K}$ with respect to 
some arbitrary ample line bundle $L$ on $X_{\bar K}$. See \cite[sec. 3]{Rossler-SSS} for the definition of Frobenius 
stabilised slopes and further references. Furthermore, the computation in 
\cite[sec. 5]{Rossler-SSS} shows that 
if 
$$
\charac(\bar K)>\rk(\Omega_{X_{\bar K}/\bar K})^2\deg_L(\Omega_{X_{\bar K}/\bar K})
$$
then we have $\bar\mu_\max(\Omega_{X_{\bar K}/\bar K})<\charac(\bar K)\cdot \bar\mu_\min(\Omega_{X_{\bar K}/\bar K}).$ Here 
$$
\deg_L(\Omega_{X_{\bar K}/\bar K}):=\int_{X_{\bar K}}{\rm c}_1(\Omega_{X_{\bar K}/\bar K})
{\rm c}_1(L)^{\dim(X_{\bar K})-1}.
$$
In our situation, we may choose 
$L=\det(\Omega_{X_{\bar K}/\bar K})$, since $\Omega_{X_{\bar K}/\bar K}$ is ample and thus 
we see that the inequality
$$
\charac(k)=\charac(\bar K)>\rk(\Omega_{X_{\bar K}/\bar K})^2\int_{X_{\bar K}}{\rm c}_1(
\Omega_{X_{\bar K}/\bar K})^{\dim(X_{\bar K})}=\rk(\Omega_{X/K})^2\int_{X}{\rm c}_1(\Omega_{X/K})^{\dim(X)}
$$
ensures that 
\begin{eqnarray*}
m_0&=&\sup\{m\in\mN|H^0(X,F_{X}^{m,*}(\Omega_{X/K}^\vee)\otimes\Omega_{X/K})\not=0\}\\
&=&
\sup\{m\in\mN|H^0(X_{\bar K},F_{X_{\bar K}}^{m,*}(\Omega_{X_{\bar K}/\bar K}^\vee)\otimes\Omega_{X_{\bar K}/\bar K})\not=0\}=0
\end{eqnarray*}
This completes Step V.

{\bf N.B.} At the time of submission of the present article, the article \cite{Rossler-SSS} by the second author (quoted in Step V) has not been accepted 
for publication yet. This article is only used in Step V and thus only the degree estimate in Theorem \ref{Mtheor} depends on it.

\end{document}